\documentclass{article}

\usepackage{latexsym}
\usepackage{epsfig}
\usepackage{amssymb}
\usepackage{amsthm}
\usepackage{amsmath}
\usepackage{epstopdf}

\def\1{{\bf 1}}

\newtheorem{thm}{Theorem}[section]
\newtheorem{cor}[thm]{Corollary}
\newtheorem{lem}[thm]{Lemma}
\newtheorem{prop}[thm]{Proposition}

\newenvironment{prf}{\noindent{\bf Proof:} }{\hfill$\Box$\mbox{}}

\title{On the spectral characterization of pineapple graphs}

\author{
Hatice Topcu\thanks{\small Dept. of Mathematics,
Nev\c{s}ehir Hac{\i} Bekta\c{s} Veli University, Turkey}\\
\small{\tt{haticekamittopcu@gmail.com}}
\\[5pt]
Sezer Sorgun\footnotemark[1]\\
\small{\tt{srgnrzs@gmail.com}}
\\[5pt]
Willem H. Haemers\thanks{\small Dept. of Econometrics and O.R.,
Tilburg University, The Netherlands}\\
\small{\tt{haemers@uvt.nl}}
}
\date{}

\begin{document}

\maketitle

\begin{abstract}

\noindent
The pineapple graph $K_p^q$ is obtained by appending $q$ pendant
edges to a vertex of a complete graph
$K_{p}$ ($q\geq 1,\ p\geq 3$).
Zhang and Zhang [{\em Some graphs determined by their spectra,}
Linear Algebra and its Applications, 431 (2009) 1443-1454]
claim that the pineapple graphs are determined by their adjacency spectrum.
We show that their claim is false by constructing graphs which
are cospectral and non-isomorphic
with $K_p^q$ for every $p\geq 4$ and various values of $q$.
In addition we prove that the claim is true if $q=2$,
and refer to the literature for $q=1$, $p=3$, and $(p,q)=(4,3)$.
\end{abstract}
AMS subject classification 05C50.

\section{Introduction}

The pineapple graph $K_p^q$ is the coalescence of the complete
graph $K_p$ (at any vertex) with the star $K_{1,q}$ at the vertex
of degree $q$.
Thus $K_p^q$ can be obtained from $K_p$ by appending $q$ pendant
edges to a vertex of $K_p$.
Clearly $K_p^q$ has $n=p+q$ vertices, ${p\choose 2}+q$ edges
and ${p\choose 3}$ triangles.
In order to exclude the complete graphs and the stars,
we assume $p\geq 3$ and $q\geq 1$.
See Figure~1 for a drawing of $K_4^4$.
Alternatively, $K_p^q$ can be defined by its adjacency matrix
\[
A=\left[\begin{array}{ccc}
0  & \1^\top   & \1^\top \\
\1 & J_{p-1}-I &   O    \\
\1 &    O      &   O
\end{array}
\right],
\]
where $\1$ is the all-ones vector (of appropriate size),
and $J_\ell$ denotes the $\ell\times\ell$ all-ones matrix.

\begin{prop}\label{cp}
The characteristic polynomial $p(x)=\det(xI-A)$ of the pineapple
graph $K_p^q$ equals
\[
p(x)=x^{q-1}(x+1)^{p-2}(x^3-x^2(p-2)-x(p+q-1)+q(p-2)).
\]
\end{prop}

\begin{prf}
The adjacency matrix $A$ has $q$ identical rows, so rank$(A)$ is
at most $p+1$, and therefore $p(x)$ has a factor $x^{q-1}$.
Similarly, $A+I$ has $p-1$ identical rows and so $(x+1)^{p-2}$
is another factor of $p(x)$.
The given partition of $A$ is equitable with quotient matrix
\[
Q=\left[\begin{array}{ccc}
0 & p-1 & q \\
1 & p-2 & 0 \\
1 &  0  & 0
\end{array}
\right]
\]
(this means that each block of $A$ has constant row sums, which
are equal to the corresponding entry of $Q$).
The characteristic polynomial of $Q$ equals
$q(x)=\det(xI-Q)=x^3-x^2(p-2)-x(p+q-1)+q(p-2)$, and
it is well known (see for example \cite{cds}, or \cite{bh})
that $q(x)$ is a divisor of $p(x)$.
\end{prf}\\

In this paper we deal with the question whether $K_p^q$ is the
only graph with characteristic polynomial $p(x)$.
In other words, is $K_p^q$ determined by its spectrum?
Note that the complete graph $K_n$ is determined by its spectrum,
but for the star $K_{1,n-1}$ this is only the case when $n=2$,
or $n-1$ is a prime.
In \cite{zz} it is stated that every pineapple graph is determined
by its spectrum.
The presented proof, however, is incorrect.
Even worse, the result is false.
In the next section we shall construct graphs with the same
spectrum as $K_p^q$ for every $p\geq 4$ and several values of $q$.

When $q=1$, the pineapple graph $K_p^1$ can be obtained from
the complete graph $K_{p+1}$ by deleting the
edges of the complete bipartite graph $K_{1,p-1}$.
Graphs constructed in this way are known to be determined by
their spectra, see~\cite{ch}.

Zhang and Zhang~\cite{zz} proved (correctly this time) that the
graph obtained by adding $q$ pendant edges to a vertex of an odd
circuit is determined by the spectrum of the adjacency matrix.
When the odd circuit is a triangle we obtain that $K_3^q$ is determined
by its spectrum.

Godsil and McKay~\cite{gm} generated by computer all pairs of
non-isomorphic cospectral graphs with seven vertices.
Since $K_4^3$ is not in their list, it is determined by its spectrum.

In Section~3 we prove that the spectrum determines $K_p^q$ when $q = 2$.
The proof uses the classification of graphs with least eigenvalue greater
than $-2$.
%This approach does not generalize.
%Because of this, and because of the existence of counter examples,
%we expect that it will be hard to settle the
%spectral characterization of the pineapple graph in general.
\\
Graphs with the same spectrum are called {\em cospectral}.
An important property of a pair of cospectral graphs is that they
have the same number of closed walks of any given length.
In particular, they have the same number of vertices,
edges and triangles.

For a graph $G$ with $n$ vertices the eigenvalues are denoted by
$\lambda_1(G)\geq\dots\geq\lambda_n(G)$.
If $H$ is an induced subgraph of $G$ with $m$ vertices, then
the eigenvalue of $H$ interlace those of $G$, which means that
$\lambda_i(G)\geq\lambda_i(H)\geq\lambda_{n-m+i}(G)$
for $i=1,\ldots,m$.
For these and other results on graph spectra we refer to \cite{cds} or \cite{bh}.

\section{Graphs cospectral with the pineapple graphs}

%Graphs with the same spectrum are called cospectral.

\begin{prop}\label{cospectral1}
Let $G$ be the graph of order $3k$ ($k\geq 2$) with adjacency matrix
\[
B=\left[
\begin{array}{c c c}
    O & J_{k}   & J_{k} \\
J_{k} & J_{k}-I & O      \\
J_{k} &    O    & J_{k}-I \\
\end{array}
\right].
\]
Then the characteristic polynomial of $G$ is equal to
\[
x^{k-1}(x+1)^{2k-2}(x-k+1)(x^2-(k-1)x-2k^2).
\]
\end{prop}

\begin{prf}
Similar to the proof of Proposition~\ref{cp}, we find that $B$ has
at least $k-1$ times an eigenvalue $0$, and at least $2k-2$ times
an eigenvalue $-1$.
The remaining three eigenvalues are the roots of the characteristic
polynomial $q(x)$ of the quotient matrix with respect to the
given partition.
It is straightforward that $q(x)=(x-k+1)(x^2-(k-1)x-2k^2)$.
\end{prf}
\\

By Proposition~\ref{cp}, we have that the characteristic polynomial
of $K_{2k}^{k^2}$ equals
\[
x^{k^2-1}(x+1)^{2k-2}(x-k+1)(x^2-(k-1)x-2k^2).
\]
So, if we add $k(k-1)$ isolated vertices to $G$ we obtain a graph
with the same characteristic polynomial, and therefore with the
same spectrum as $K_{2k}^{k^2}$.
For $k=2$, the two cospectral graphs are drawn in Figure~1.
\begin{figure}[h]
\hspace{20pt}
\includegraphics[height=160pt,width=350pt]{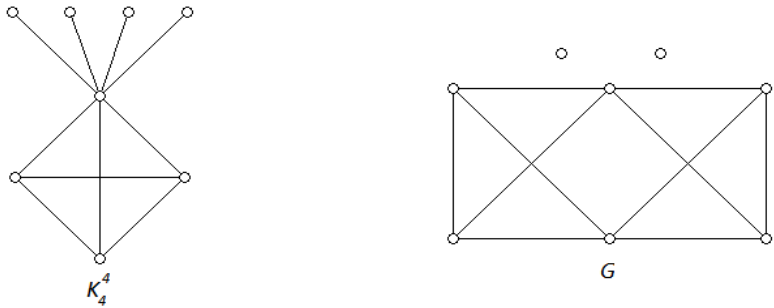}
\caption{Two graphs with characteristic polynomial
$x^3(x+1)^2(x+2)(x^2-x-8)$}
\end{figure}%

%\begin{cor}
%If $p=2k\geq 4$ and $q=k^2$, then $K_p^q$ is not determined by
%the spectrum of its adjacency matrix.
%\end{cor}
%

In the next result we obtain more graphs which are nonisomorphic but cospectral
with the pineapple graph.
It makes use of the graph $K_n\!\setminus\! K_m$ ($m<n$), which
can be obtained from the complete graph $K_n$ be deleting the edges
of a subgraph $K_m$ (in other words, $K_n\!\setminus\! K_m$ is the
complete multipartite graph $K_{m,1,\ldots,1}$ with $n$ vertices).
The disjoint union of two graphs $G$ and $H$ is denoted by $G+H$.
\begin{prop}\label{cospectral2}
Let $k\geq 2$ and $p$ be integers such that $r=k(k-1)/(p-k-1)$ is a positive integer.
Then
\[
K_{p+r}\!\setminus\! K_{k+r} + K_k + k(k-2)K_1\ \mbox{ and }\ K_p^{r(p-k)}
\]
are cospectral.
The common characteristic polynomial equals
\[
x^{r(p-k)-1}(x+1)^{p-2}(x-k+1)(x^2-x(p-k-1)-(k+r)(p-k)).
\]
\end{prop}
\begin{prf}
It is well known (see for example~\cite{ch}) and easily proved that
the characteristic polynomial of $K_n\!\setminus\! K_m$ is equal to
$x^{m-1}(x+1)^{n-m-1}(x^2-(n-m-1)x-m(n-m))$.
From this and Proposition~\ref{cp} we find that $K_{p+r}\!\setminus\! K_{k+r} + K_k + k(k-2)K_1$
and $K_p^{r(p-k)}$ both have the characteristic polynomial given above.
\end{prf}
\\

For fixed $k\geq 2$ the above condition is fulfilled for several values of $p$.
Every $k\geq 6$, gives at least eight values of $p$ for which the condition is satisfied.
In particular it follows that $K_p^q$ is not determined by its spectrum when
$q=2(p-2)(p-3)$, $p\geq 4$, and when $q=3(p-3)(p-4)/2$, $p\geq 5$ (including $K_5^3$).
%In particular, we find that $K_5^3$ is not determined by its spectrum.
%In particular, $1$, $2$, $\lfloor k/2 \rfloor$, $k-1$, $k$, $k(k-1)/\lfloor k/2 \rfloor$, $k(k-1)/2$, and $k(k-1)$ are
%divisors of $k(k-1)$, which are all distinct when $k\geq 6$.
Moreover, $p=2k\geq 4$, $q=k^2$ satisfies the requirements of both propositions above.
Thus we have:
\begin{cor}
If $p$ is even, $p\geq 4$ and $q=(p/2)^2$, then there exist at least three nonisomorphic graphs with the spectrum of $K_p^q$.
\end{cor}

\section{The spectral characterization of $K_p^2$}

The main tool in this section is the classification of graphs
with least eigenvalue strictly greater than $-2$, due to Doob
and Cvetkovi\'c~\cite{dc}.
Below we present this classification.
We inserted results from \cite{cl} that relate the different cases
with the so-called discriminant $d_G$ of a graph $G$, defined by
$d_G=|p(-2)|$, where $p(x)$ is the characteristic polynomial of $G$.

The classification uses the not so well-known notion of a
generalized line graph of type $(1,0,\ldots,0)$,
so we will recall it here.
Suppose $G$ is a graph, then the {\em line graph} $L(G)$ is the
graph whose vertex set consists of the edges of $G$, where two
edges are adjacent in $L(G)$ whenever they intersect in $G$.
Assume the vertices of $G$ are labeled $1,2,\ldots,n$, then the
{\em generalized line graph of $G$ of type $(1,0,\ldots,0)$}
consists of $L(G)$ extended with two non-adjacent new vertices
which are adjacent to all vertices of $L(G)$ that correspond to
an edge of $G$ incident with the vertex with label~$1$.
For example, $K_p^2$ is the generalized line graph of $K_{1,p}$ of
type $(1,0,\ldots,0)$, provided a vertex of degree $1$ has label~$1$.
\begin{lem}\label{cd}
Let $G$ be a connected graph with least eigenvalue greater than $-2$.
Then one of the following statements holds:
\begin{itemize}
\item[(i)] $G$ has eight vertices, and $d_{G}=1$.
\item[(ii)] $G$ has seven vertices, and $d_{G}=2$.
\item[(iii)] $G$ has six vertices, and $d_{G}=3$.
\item[(iv)] $G$ is the line graph of an unicyclic graph with a
cycle of odd length, and $d_{G}=4$.
\item[(v)] $G$ is a generalized line graph of a tree of type
$(1,0,\ldots,0)$, and $d_G=4$.
\item[(vi)] $G$ is the line graph of a tree with $n\geq 5$
vertices, and $d_{G}=n$.
\end{itemize}
\end{lem}

\begin{thm}\label{q=2}
The pineapple graph $K_{p}^{2}$ is determined by the spectrum of
its adjacency matrix.
\end{thm}

\begin{prf}
Suppose that $G$ is a graph cospectral with $K_p^2$.
Then $G$ has $p+2$ vertices, ${p\choose 2}+2$ edges and
$p\choose 3$ triangles,
and the characteristic polynomial $p(x)$ of $G$ is equal to
\[
x(x+1)^{p-2}(x^{3}-(p-2)x^{2}-(p+1)x+2(p-2)).
\]
From this it follows that $\lambda_{2}(G)<1$,
and $\lambda_{p+2}(G)>-2$.

Suppose $G$ is disconnected.
If there exist two components with at least one edge, then $2K_2$ is
an induced subgraph of $G$.
Eigenvalue interlacing gives $1=\lambda_2(2K_2)<\lambda_2(G)<1$,
contradiction.
So only one component of $G$ has edges.
The multiplicity of the eigenvalue $0$ equals $1$, therefore $G$
has two components, one of which is an isolated vertices.
The larger component $G'$ has characteristic polynomial $p(x)/x$,
so $d_{G'}=|p(-2)/2|=2$.
Hence, by Lemma~\ref{cd}, $G'$ has seven vertices.
We know that $G'$ has $17$ edges and $20$ triangles.
One straightforwardly checks that there are exactly ten
connected graphs on seven vertices with 17 edges,
and none has 20 triangles.
This proves that $G$ is connected and $d_G=p(-2)=4$.
So we are dealing with case $(iv)$ or $(v)$ of Lemma~\ref{cd}.
\\

In case~$(iv)$, $G$ is the line graph of a unicyclic graph $H$
with $p+2$ vertices.
Let $v$ be a vertex of $H$ of maximum degree $k$ and put $\ell=p+2-k$.
When $v$ is not in a triangle of $H$, every edge not incident with
$v$ is incident with at most one neighbor of $v$.
Therefore there are at least $\ell(k-1)$ pairs of disjoint edges
in $H$.
If $v$ is in a triangle $\{u,v,w\}$ then the edge $\{u,w\}$ is incident with
$u$ and $w$ and disjoint from $k-2$ edges through $v$.
So we obtain at least $\ell(k-1)-1$ pairs of disjoint edges in $H$.
Disjoint edges in $H$ correspond to pairs of nonadjacent vertices in $G$.
This implies that $G$ has at most ${{p+2}\choose 2} - \ell(k-1)+1$ edges.
Thus we find ${{p} \choose 2}+2 \leq {{p+2}\choose 2}-\ell(k-1)+1$,
which leads to $(\ell-2)(k-3)\leq 2$.
This implies that $k\leq3$, $k\geq p$, or
$(p,k)\in\{(5,4),(6,4),(6,5)\}$.

If $k\leq 3$, then $H$ has at most $(p+2)/2$ vertices of degree~$3$
(indeed, the average degree equals $2$, so the number of vertices of degree 1
equals the number of vertices of degree 3).
Therefore $G$ has at most $1+(p+2)/2$ triangles,
and therefore ${p\choose 3}\leq 1+(p+2)/2$, which gives $p\leq 4$.
The case $p=3$ is solved in \cite{zz}, and when $p=4$,
$G$ has $8$ edges and $4$ triangles and no $K_4$ (since $H$ has
maximum degree $3$), which is impossible.

If $k=p+1$, then $H=K_3^{p-1}$, hence $G$ has
${{p+1}\choose 2}+2$ edges, which is false.

If $k=p$, then there are just two ways to create an odd cycle in $H$
by adding two edges.
In both cases we obtain $K_3^{p-2}$ with one pendent edge attached to
a vertex of degree at most $2$ which leads to more than
${p}\choose 3$ triangles in $G$.

Suppose $(p,k)\in\{(5,4),(6,4),(6,5)\}$.
We get the maximum possible number of triangles in $G$ if $H$
consist of a triangle with $k-2$ pendent edges in one vertex of
the triangle and $p-k+1$ pendant edges on another vertex.
In all three cases $G$ has fewer triangles than $K_{p}^2$.

So we can conclude that there are no graphs cospectral with $K_{p}^2$
in case $(iv)$.
\\

Finally, we consider case~$(v)$ of Lemma~\ref{cd}.
Then $G$ is a generalized line graph of a tree $T$ with $p+1$
vertices and $p$ edges of type $(1,0,\ldots,0)$.
Let $H$ be the line graph of $T$, let $v$ be a vertex of $T$ with maximum degree $k$,
and define $\ell=p+2-k$.
Then $H$ has at most ${k\choose 2}+{{\ell-1}\choose 2}$ edges,
and $G$ has at most ${k\choose 2}+{{\ell-1} \choose 2}+2k$ edges,
hence
${k\choose 2}+{{\ell-1} \choose 2}+2k \geq {{p}\choose 2}+2$,
which leads to $(\ell-4)(k-1)\leq 0$.
So $p=k+2$, $p=k+1$, or $p=k$.

If $p=k+2$,then equality holds in the above inequality, which
means the two edges of $T$ which are not incident with $v$ intersect.
This implies that $H$ is the coalescence of $K_k$ and a triangle,
or the coalescence of $K_k$ with the path $P_3$ at a pendent vertex of $P_3$.
This leads to seven possible generalized line graphs of type $(1,0,\ldots,0)$
(depending on which vertex has label~$1$).
None of these has $p\choose 3$ triangles.

If $p=k+1$, then $H$ consists of the complete graph $K_k$ with one pendent edge.
Now there are only three possible generalized line graphs of the required type,
and none has $p\choose 3$ triangles.

If $p=k$, then $T=K_{1,{k}}$, and $H=K_k$.
There are two generalized line graphs possible,
but only $K_p^2$ has $p\choose 3$ triangles.
So $G=K_p^2$.
\end{prf}

\section{Concluding remarks}
For $q\geq 3$, the smallest eigenvalue of $K_p^q$ is less than $-2$,
so the proof of Theorem~\ref{q=2} does not generalize to larger $q$.
Moreover, since $K_4^4$ and $K_5^3$ are not determined by their spectrum,
it will already be difficult to settle the spectral characterization of $K_p^q$
when $p=4$ and when $q=3$.

The graphs cospectral with $K_p^q$ presented in Propositions~\ref{cospectral1}
and~\ref{cospectral2} have an integral eigenvalue $k-1$.
So it is conceivable that the pineapple graphs with three nonintegral eigenvalues are
determined by their spectrum.
For such pineapple graphs the three nonintegral eigenvalues belong to one component,
so the graph consists of one large component and a number of isolated vertices.
But we doubt if this observation is useful, since the examples of Proposition~\ref{cospectral1}
also have this structure.

This work is supported by TUBITAK (the scientific and technological research council of Turkey) 2214-A doctorate research fellowship program.

%We believe that for a complete solution of the spectral characterization of the pineapple graph
%(if at all possible) a new approach is needed.
%Perhaps it is possible to classify
%all graphs with all but three eigenvalues equal to $0$ or $- 1$,
%but that would be looks like a long shot.

\end{document}